\documentclass[conference]{IEEEtran}
\IEEEoverridecommandlockouts
\usepackage[table]{xcolor}
\usepackage{cite}
\usepackage{amsmath,amssymb,amsfonts,bm,bbm}
\usepackage{graphicx}
\usepackage{textcomp}
\usepackage{xcolor}
\usepackage[stable]{footmisc}

\def\BibTeX{{\rm B\kern-.05em{\sc i\kern-.025em b}\kern-.08em
    T\kern-.1667em\lower.7ex\hbox{E}\kern-.125emX}}
    
\usepackage[utf8]{inputenc}
\usepackage{url}

\usepackage{amsmath}
\usepackage{mathtools,amssymb,lipsum,nccmath,bm,amsmath}

\usepackage{cuted}
\usepackage{array}
\usepackage{ragged2e}
\usepackage{optidef}
\usepackage{algorithm,algcompatible}
\usepackage{siunitx}

\algnewcommand\INPUT{\item[\mathbf{Inputs:}]}
\algnewcommand\OUTPUT{\item[\mathbf{Output:}]}
    
\usepackage{dblfloatfix}
\usepackage{float}
  \floatstyle{ruled}
  \newfloat{algorithm}{htbp}{loa}
  \floatname{algorithm}{Algorithm}
\usepackage{url}

\usepackage{balance}%
\usepackage{amsthm}

\begin{document}

\title{Power System Quasi-Steady State Estimation: \\ An Echo State Network Approach\footnote{We thank Professor Douglas Plaza at ESPOL Polytechnic University for his assistance with the traditional methods.}}

\author{
\IEEEauthorblockN{Gabriel Intriago\IEEEauthorrefmark{1}, Holger Cevallos\IEEEauthorrefmark{2}, and Yu Zhang\IEEEauthorrefmark{1}}
%\vspace{0.05in}
\IEEEauthorblockA{\IEEEauthorrefmark{1}University of California, Santa Cruz. Emails: \emph{\{gintriag, zhangy\}@ucsc.edu}}
%\vspace{0.05in}
\IEEEauthorblockA{\IEEEauthorrefmark{2}ESPOL Polytechnic University. Email: \emph{hcevallo@espol.edu.ec}}
\thanks{This work was supported in part by a Seed Fund Award from CITRIS and the Banatao Institute at the University of California, and the Hellman Fellowship. We thank Professor Douglas Plaza at ESPOL Polytechnic University for his assistance with the traditional methods.}
%\vspace{0.05in}
}

\maketitle

\begin{abstract}
The operating point of a power system may change due to slow enough variations of the power injections. Rotating machines in the bulk system can absorb smooth changes in the dynamic states of the system. In this context, we present a novel reservoir computing (RC) method for estimating power system quasi-steady states. By exploiting the behavior of an RC-based recurrent neural network, the proposed method can capture the inherent nonlinearities in the power flow equations. Our approach is compared with traditional methods, including least squares, Kalman filtering, and particle filtering. We demonstrate the estimation performance for all the methods under normal operation and sudden load change. Extensive experiments tested on the standard IEEE 14-bus and 300-bus cases corroborate the merit of the proposed approach.
\end{abstract}


\section{Introduction}
This paper introduces a framework for power systems quasi-steady-state estimation (QSSE) based on reservoir computing (RC). RC is a scheme derived from recurrent neural networks (RNN) that maps the inputs of a system into a high-dimensional space called reservoir \cite{Ghosh2019}. The echo state network (ESN) has a fixed, random, and sparsely connected reservoir, which is an RNN that belongs to the RC framework. Different from traditional RNNs, the ESNs are easier to train with guaranteed converge if the echo state property holds \cite{Jaeger2001}.

A quasi-steady state operation condition refers to the system operating point driven by slow and gradual load or renewable generation changes. In other words, the changes in the dynamic states are negligible due to the generators absorbing the changes instantly; such a situation makes it interesting to track the changes of algebraic quantities, i.e., the voltage phasors \cite{Zhao2019}. The QSSE approach provides acceptable results when the trajectory of the state progresses smoothly. Under certain conditions, the mathematical formulation of the state transition model is assumed to be linear and derived from time series forecasting techniques such as exponential smoothing \cite{LeitedaSilva1983}. However, QSSE exhibits inaccurate estimations in the presence of abrupt load changes because it may take some time for the parameters of the linearized state transition model to adjust to the new scenario. We hypothesize that the ESN can adapt to such a dynamic scenario while catching the nonlinear relationship from the power flow equations.

% The ESNs have proven outstanding results in neuroscience \cite{Grel2009}, signal processing \cite{Kudithipudi2016}, and wireless communications \cite{Jaeger2004}. 

In general, solutions based on artificial neural networks have penetrated several fields such as navigation \cite{Pico_mdpi}, robotics \cite{Pico2023}, estimation \cite{Pico2022}, telecommunications \cite{Bakhadirovna2021}, and control \cite{Fei2022}. The ESN, a type of recurrent neural network, has proven outstanding results in neuroscience \cite{Grel2009}, signal processing \cite{Kudithipudi2016}, and wireless communications \cite{Jaeger2004}. Mainly, ESNs have been considered for several power system applications such as true harmonic current detection \cite{Dai20092}, wide-area monitoring \cite{Venayagamoorthy2007}, and power system nonlinear load modeling \cite{JingDai2008}. Broadly speaking, \cite{Goswami2021} explores the performance of an ensemble Kalman filter (KF) in the feedback path to the reservoir of an ESN by using sparse measurements for a complex nonlinear system. An ESN combined with an extended KF is developed for real-time modeling prediction for ship motion \cite{Peng2017}. An online prognostic method, which combines an ESN with a kernel recursive least squares and a Bayesian technique, is developed for tracking the health status of a degraded system \cite{Zhou2018}. The contributions of this paper are two-fold:
\begin{itemize}
    \item We develop a novel reservoir-computing-based state estimation method for power systems under the assumption that the changes in dynamical states are negligible. To the best of the authors' knowledge, this is the first time a reservoir computing network has been proposed for power systems state estimation.
    \item Our proposed approach is comparable in time and accuracy to the industry standard, the Weighted Least Squares, and other more complex filtering methods such as the extended Kalman Filter and the particle filter.
\end{itemize}

% This document is organized as follows. Section II presents the framework for state estimation and the different estimation methods. Section III shows the simulations and numerical results. Finally, section IV discusses the conclusion ideas.

\section{{State estimation methods}}\label{sec:methods}

\subsection{Framework for QSSE} % \cite{Venkata2011}
Generally, the system model for state estimation is characterized as a continuous-time state-space model. However, the model is transformed into its discrete-time state-space form in practice through a discretization technique. More formally, a general discrete-time state-space model for quasi-steady-state estimation is presented as follows \cite{Zhao2019}
\begin{subequations}
\begin{align}
\bm{x}_k &= \bm{f}\left(\bm{x}_{k-1},\bm{u}_k,\bm{p}\right) + \bm{w}_k,\label{eq1}\\
\bm{z}_k &= \bm{h}\left(\bm{x}_k,\bm{u}_k,\bm{p}\right) + \bm{v}_k\label{eq2},
\end{align}
\end{subequations}
where $k$ is the discrete-time step index; $\bm{x}_k\in\mathbb{R}^n$ represents the voltage phasor state vector; $\bm{u}_k\in\mathbb{R}^u$ is the system input vector; $\bm{p}\in\mathbb{R}^p$ represents the model parameters; $\bm{z}_k\in\mathbb{R}^m$ is the measurement vector; $\bm{f}$ is the nonlinear process function; $\bm{h}$ is the nonlinear measurement function; $\bm{w_k}$ and $\bm{v_k}$ are often assumed to be uncorrelated and normally distributed with zero mean and covariance matrices $\bm{Q}_k$ and $\bm{R}_k$, respectively. The time index $k$ for both covariance matrices is dropped; i.e., $\bm{Q}$ and $\bm{R}$ remain fixed at all time steps. In this work, we omit the vectors $\bm{u}_k$ and $\bm{p}$, given that we do not control the system, and the parameters are assumed to be time-invariant.

% WLS
\subsection{{The Weighted Least Squares (WLS)}} \label{sec:wls}
The WLS estimator of PSSE minimizes the following objective at each time step \cite{Abur2004} \cite{Cevallos20181}
\begin{align}
 J(\bm{x}_k)
 =[\bm{z}_k-\bm{h}(\bm{x}_k)]^{\top}\bm{R}^{-1}[\bm{z}_k-\bm{h}(\bm{x}_k)].
\end{align}
The first-order optimality conditions must be satisfied at the minimum, which yields
\begin{equation}
 \bm{g}(\bm{x}_k)=\frac{\partial J(\bm{x}_k)}{\partial \bm{x}_k}=-\bm{H}^\top(\bm{x}_k) \bm{R}^{-1}[\bm{z}_k-\bm{h}(\bm{x}_k)]=\bm{0},
\end{equation}
where $ \bm{H}(\bm{x}_k)=\frac{\partial \bm{h}(\bm{x}_k)}{\partial \bm{x}_k}
$ is the Jacobian matrix of $\bm{h}(\bm{x}_k)$.

We obtain the Gauss-Newton update by expanding the non-linear function $\bm{g}(\bm{x})$ into its first-order Taylor expansion around the state $\bm{x}_k$
% \begin{equation}
%     \bm{g}(\bm{x})=\bm{g}(\bm{x}_k)+\bm{G}(\bm{x}_k)(\bm{x}-\bm{x}_k)=0.
% \end{equation}
\begin{equation}
 \bm{x}_{k}^{i+1}=\bm{x}_k^i-\bm{G}^{-1}(\bm{x}_k^i)\bm{g}(\bm{x}_k^i),
\end{equation}
where 
 $\bm{x}_k^i$ is the solution at iteration $i$ and
 \begin{subequations}
 \begin{align}
 \bm{G}(\bm{x}_k^i)&=\frac{\partial \bm{g}(\bm{x}_k^i)}{\partial \bm{x}_k}=\bm{H}^\top(\bm{x}_k^i)\bm{R}^{-1}\bm{H}(\bm{x}_k^i), \\ \bm{g}(\bm{x}_k^i)&=-\bm{H}^\top(\bm{x}_k^i)\bm{R}^{-1}[\bm{z}_k-\bm{h}(\bm{x}_k)].
 \end{align}
 \end{subequations}
To bypass the matrix inversion, the incremental update $\Delta\bm{x}_{k}^{i+1} = \bm{x}_{k}^{i+1} - \bm{x}_{k}^{i}$ is often obtained by solving the normal equation
$
 \bm{G}(\bm{x}_k^i)\Delta\bm{x}_{k}^{i+1} = \bm{H}^\top(\bm{x}_k^i)\bm{R}^{-1}[\bm{z}_k-\bm{h}(\bm{x}_k)]
$.

% EKF
\subsection{{The Extended Kalman Filter (EKF)}} \label{sec:ekf}
The Kalman filter framework consists of two iterative steps: a prediction step based on \eqref{eq1} and a filtering/update step based on \eqref{eq2}. Typically, given the state estimate $\bm{x}_{k-1}$ and its covariance matrix $\bm{P}_{k-1}$, the predicted state $\bm{x}_k^{-}$ is obtained from \eqref{eq1} directly or through a set of points drawn from the probability distribution of $\bm{x}_{k-1}$, which depends on the distribution of $\bm{w}_k$. The prediction state is combined with measurements to estimate the state vector and its covariance matrix for the update step.

The EKF linearizes the system's model around the current operating point by using the Jacobian matrices of $\bm{f}$ and $\bm{h}$ \cite{Cevallos2018UTE}. We assume that any change in the operating point is driven by slow enough stochastic changes in the power injections. The dynamics of the synchronous machines and loads are sufficiently small and can be neglected. Such an assumption allows us to drop the need for the Jacobian matrix of $\bm{f}$ and introduce the forecasting-aided version of \eqref{eq1} in the EKF prediction step as follows
\begin{itemize}
    \item \emph{EKF prediction step}:
    \begin{subequations}
    \begin{align}
        \bm{x}_k^{-} &= \bm{F}_{k-1}\bm{x}_{k-1}+\bm{g}_{k-1}+\bm{w}_{k-1}\\
        \bm{P}_{k}^{-} &= \bm{F}_{k-1}\bm{P}_{k-1}\bm{F}_{k-1}^\top + \bm{Q}.
    \end{align}
    \end{subequations}
\end{itemize}
$\bm{F}$ and $\bm{g}$ are computed online based on Holt's linear exponential smoothing technique \cite{LeitedaSilva1983}
\begin{subequations}
\begin{align}
    \bm{F}_k &= \alpha_k(1+\beta_k)\mathbf{I},\\
    \bm{g}_k &= (1+\beta_k)(1-\alpha_k)\bm{x}_k^{-}- \notag \\
    &\hspace{0.5cm}\beta_k\bm{a}_{k-1} + (1-\beta_k)\bm{b}_{k-1}, \\
      \bm{a}_k &= \alpha_k\bm{x}_k+(1-\alpha_k)\bm{x}_k^{-},\\
    \bm{b}_k &= \beta_k(\bm{a}_k-\bm{a}_{k-1})+(1-\beta_k)\bm{b}_{k-1},
\end{align}
\end{subequations}
where coefficients $\alpha_k, \beta_k \in [0,1]$.

\begin{itemize}
    \item \emph{EKF update step}:
    \begin{subequations}
    \begin{align}
        \bm{x}_k &= \bm{x}_k^{-} + \bm{K}_k\left(\bm{z}_k - \bm{h}(\bm{x}_k)\right),\\
        \bm{K}_k &= \bm{P}_k^{-}\bm{H}_k^\top\left(\bm{H}_k\bm{P}_k^{-}\bm{H}_k^\top+\bm{R}\right)^{-1},\\
        \bm{P}_k &= \left(\bm{I}-\bm{K}_k\bm{H}_k\right)\bm{P}_k^{-}.
    \end{align}
    \end{subequations}
\end{itemize}
$\bm{x}_k^{-}$ and $\bm{P}_k^{-}$ are known as \textit{a priori} mean and covariance matrix of the state. $\bm{x}_k$ and $\bm{P}_k$ are the \textit{a posteriori} counterparts. The Kalman gain is the matrix $\bm{K}_k$.

\subsection{{Particle Filter (PF)}} \label{sec:pf}
Assuming that the system can be modeled as a stochastic process, the PF framework finds the best estimate for the true state $\bm{x}_k$ via a set of particles (a.k.a. samples) that represent the posterior distribution of the states given the measurements $\bm{z}_k$. The approximate posterior distribution is given as
\begin{align}
    p(\bm{x}_{k-1}|\bm{z}_{1:k-1}) &\approx \sum_{i=1}^N{W}_{k-1}^i\delta\left(\bm{x}_{k-1} - \bm{x}_{k-1}^i\right).
\end{align}

Without any assumptions about the state-space model or the state distribution, the PF is governed by the following steps \cite{Zhou2015} \cite{Intriago2018}
\begin{itemize}
    \item \emph{PF prediction step}:
        \begin{align}
            \bm{x}_k^{i-} = \bm{f}(\bm{x}_{k-1}^{i}) + \bm{w}_{k-1}^i.
        \end{align}
    \item \emph{PF update step}:
    \begin{subequations}
        \begin{align}
        \tilde{W}_k^i &= {W}_{k-1}^i\times p_{vk}\left(\bm{z}_k|\bm{h}\left(\bm{x}_k^{i-}\right)\right),\\
        {W}_{k}^i &= {\tilde{W}_k^i}\mathbin{/}{\sum_{j=1}^N\tilde{W}_k^j},\\
        \bm{x}_k &= \sum_{i=1}^N\bm{x}_k^{i-}{W}_k^i.
     \end{align}
    \end{subequations}
    \item \emph{PF resampling step}:
    Optionally at each time $k$, we take $N$ samples with replacement from the set $\left\{\bm{x}_k^{i-}\right\}_{i=1}^N$, where the probability to take sample $i$ is $W_k^i = 1/N$ \cite{Gustafsson2010}.
\end{itemize}

$\tilde{W}_k^i$ represents the updated weight of the $i$th sample while ${W}_k^i$ is the normalized updated weight. The likelihood of the measurements $\bm{z}_k$ given the noiseless measurement is represented by  $p_{vk}\left(\bm{z}_k|\bm{h}\left(\bm{x}_k^{i-}\right)\right)$. $N$ is the number of particles used to approximate the posterior distribution of the state.

% ECHO STATE NETWORKS
\subsection{Echo State Networks (ESN)}
\begin{figure}
    \centering
     \includegraphics[width=0.7\linewidth]{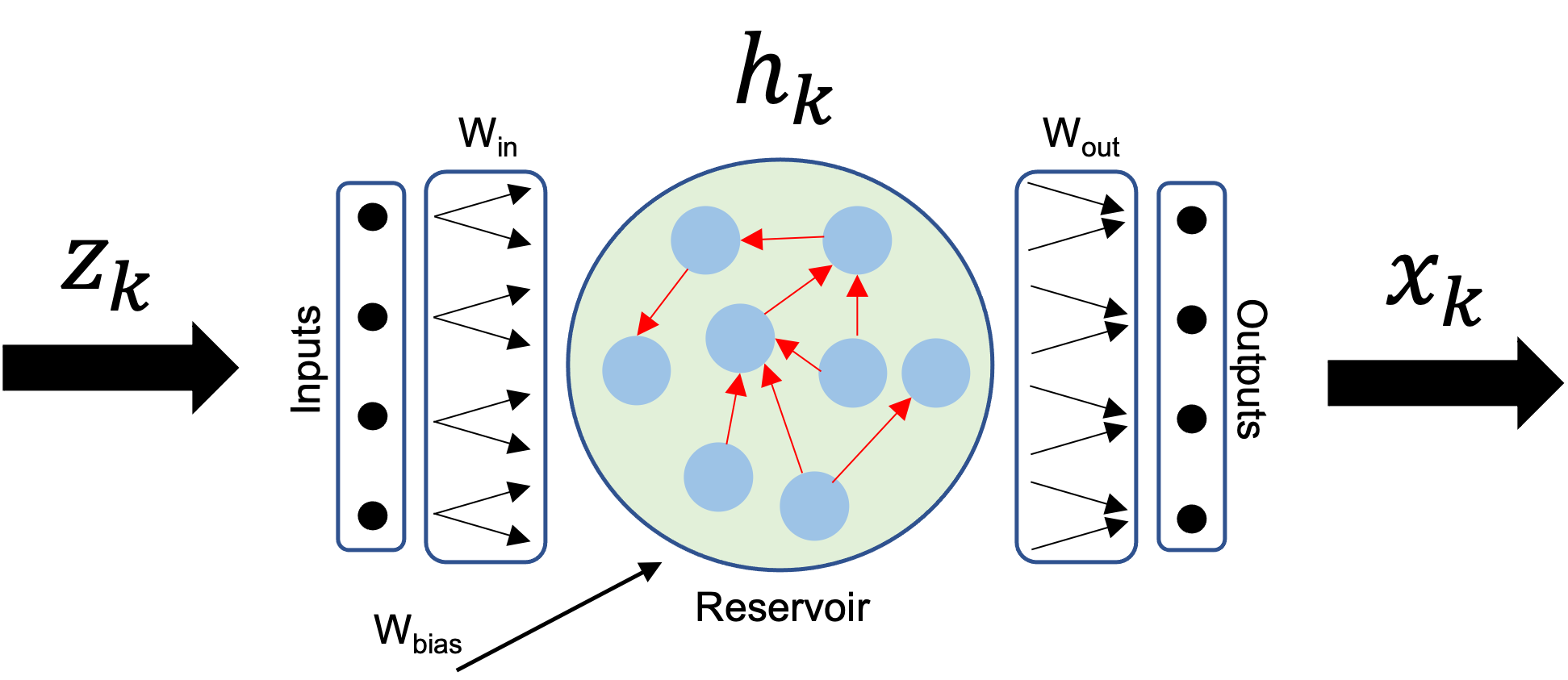}
     \caption{The architecture of an echo state network.}
     \label{fig1esn}
\end{figure}

In this paper, we implement an ESN for the QSSE of power systems, which provides a structure and supervised learning principle for recurrent neural networks (RNNs). The general idea behind the ESN is to drive a sparsely connected reservoir of neurons with the input signal and then obtain the desired output signal by combining the neurons' responses  \cite{Jaeger2001}. Fig.~\ref{fig1esn} shows the general architecture of an ESN.

For the QSSE task, the measurement vector $\bm{z}_k\in\mathbb{R}^{m}$ is the ESN input vector while the state vector $\bm{x}_k\in\mathbb{R}^{n}$ is the desired ESN output signal. Let $\bm{h}_k\in\mathbb{R}^{N_h}$ denote the reservoir's internal signal; $\bm{W}^{\text{res}}\in\mathbb{R}^{N_h\times N_h}$ the reservoir weight sparse matrix with spectral radius $\rho_{max}(\bm{W}^{\text{res}}) < 1$ and elements following a normal distribution
centered around zero; ${\bm{W}^{\text{in}}}\in\mathbb{R}^{N_h\times m}$ and ${\bm{W}^{\text{bias}}}\in\mathbb{R}^{N_h}$ the input and bias weight matrices with uniformly distributed entries; and ${\bm{W}^{\text{out}}} \in\mathbb{R}^{n\times N_s}$ is the output weight matrix, where $N_s = 1+m+N_h$. Then, the internal signal is updated as follows
\begin{equation}
\bm{h}_{k+1} = (1-\alpha)\bm{h}_{k} + \alpha \mathbf{f}_1\big(\bm{W}^{\text{res}}\bm{h}_{k} + \bm{W}^{\text{in}}\bm{z}_k + \bm{W}^{\text{bias}}\big),\label{eqesn1}
\end{equation}
where {$\mathbf{f}_1$} is an invertible nonlinear function (e.g., sigmoid or tanh). Parameter $\alpha$ controls how sensitive the reservoir is towards the past or current values. The desired output is then obtained by
\begin{align}
\bm{x}_k &= \mathbf{f}_2\left(\bm{W}_k^{\text{out}}\bm{s}_k\right),\\
\bm{s}_k &= [1, \bm{z}_k^\top, \bm{h}_k^\top]^\top\label{eqesn2}
\end{align}
where
$\bm{s}_k\in\mathbb{R}^{N_s}$ is the augmented internal signal. Typically, {$\mathbf{f}_2$} is the identity mapping.

We leverage the online version of the ESN commonly trained by two strategies: i) the incremental ridge regression and ii) the least mean squares (LMS) filter. The former technique trains the ESN by solving the following problem
\begin{align}
\min_{\bm{W}^{\text{out}}} \, \frac{1}{n}\sum_{i=1}^{n}\left(\sum_{k=1}^T\left(x_{k_i} - x_{k_i,o}\right)^2 + \epsilon\left\|\mathbf{w}_{i}^{\text{out}}\right\|^2\right).
\label{optridge}
\end{align}
The decision variable is $\bm{W}^{\text{out}}
=[{\mathbf{w}_{1}^{\text{out}}}^{\top},\dots, {\mathbf{w}_{n}^{\text{out}}}^{\top}]^{\top}$, where $\mathbf{w}^{\text{out}}_{i}$ denotes the $i$-th row; $T$ is the time horizon; $x_{k_i}$ is the $i$th component of the estimate state vector at step $k$ while $x_{k_i,o}$ is the true counterpart; and $\epsilon$ is the regularization constant.
The most universal and stable solution to problem \eqref{optridge} is given as
\begin{align}
    \bm{W}^{\text{out}} &= \bm{X}_{o}\bm{H}^\top\left(\bm{H}\bm{H}^\top + \epsilon\mathbf{I}\right)^{-1}, \label{eq:stableSolu}
\end{align}
where $\bm{X}_{o}\in\mathbb{R}^{n\times T}$ is the matrix containing all the true states vectors from $k=1$ up to $k=T$. Matrix $\bm{H}=[\bm{s}_1,\dots,\bm{s}_T]\in\mathbb{R}^{N_s\times T}$ collects all the reservoir's augmented internal signal vectors across time. Note that $\bm{X}_{o}\bm{H}^\top \in \mathbb{R}^{n\times N_s}$ and $\bm{H}\bm{H}^\top \in \mathbb{R}^{N_s\times N_s}$ do not depend on the time horizon $T$, and hence can be incrementally updated at each time step \cite{Lukosevicius2012}.

The LMS filter trains the ESN based on the gradient descent algorithm. The algorithm starts by assuming small weights (zero in most cases) and, at each time step, the weights are updated by finding the gradient of the mean square error \cite{Widrow1988} 
\begin{subequations}
        \begin{align}
        \bm{e}_k & = \bm{x}_{k,o} - \bm{x}_k \\ 
        J_k &= \frac{1}{2}\| \bm{e}_k \|^2 \\
        \bm{W}_{k+1}^{\text{out}} &= \bm{W}_{k}^{\text{out}} - \eta \nabla J_k \\
        \bm{W}_{k+1}^{\text{out}} &= \bm{W}_{k}^{\text{out}} + \eta \bm{e}_k\bm{s}_k^\top
        \label{eqlms}
        \end{align}
\end{subequations}
where $\eta$ is the learning rate of the LMS filter. 

Our proposed approach combines the incremental ridge regression technique with the LMS filter. At step $k+1$, we incrementally update the output weight matrix $\bm{W}_{k+1}^{\text{out}}$ via equation \eqref{eqlms}. The first term $\bm{W}_{k}^{\text{out}}$ is obtained from the regularized incremental ridge regression update in \eqref{eq:stableSolu}.

\begin{table}[t]
\caption{ESN parameters for voltage phasors estimation}
\begin{center}
\begin{tabular}{c||c|c} 
 \hline
 % &&&&\\
ESN&Voltage&Voltage\\
 %Ratio&(\%)&(\%)&(\%)&(\%)\\
Parameter& Magnitude & Angle \\
 %&$k = 10$&$k = 40$&$k = 80$&$k = 130$\\
 \hline
 \hline
 Reservoir Size&400&21\\
 \hline
 Bias scaling&$0.75 \times 10^{3}$&$0.6 \times 10^{3}$\\
 \hline
 Input scaling&$1 \times 10^{-5}$&$1 \times 10^{-5}$\\ 
 \hline
 Output scaling&$1 \times 10^{-5}$&$1 \times 10^{-5}$\\ 
 \hline
 $\epsilon$&$1 \times 10^{-9}$&$1.5 \times 10^{-5}$\\
 \hline
 $\alpha$&$1.63 \times 10^{-2}$&$2 \times 10^{-1}$\\ 
 \hline
 Spectral radius&0.3&0.05\\ 
 \hline
 $\eta$&$0.9 \times 10^{-2}$&$0.5 \times 10^{-2}$\\ 
 \hline
 \end{tabular}
\end{center}
\label{tbesnparam}
\vspace*{-0.5cm}
\end{table}

%------------------------- EXPERIMENTS AND RESULTS
% , which can be seen in Fig. \ref{fig2}
\section{Experiments and Results}
\subsection{Simulations Description}
The methods presented in Section~\ref{sec:methods} were tested on the IEEE 14-bus and IEEE 300-bus test cases using MATLAB and the load flow analysis library from \cite{Panda1998}. One hundred time-sample intervals were obtained by running load flows under smooth variations in the loading conditions to simulate that any change in the system's operating point is sufficiently slow. In some selected buses, the loads change with a linear trend of $1\%$ over the whole simulation time horizon. We consider measurements such as voltage magnitudes, active/reactive power injections, and active/reactive power flows. In addition, voltage magnitudes, power injections, and power flow measurements are injected with a zero-mean random additive Gaussian noise with a standard deviation $\sigma$ of $1\%$ for voltage magnitudes and $2\%$ for powers. The initialization of the ESN is carried out by using the first 30 time samples. The estimation process starts at time step $k=31$ and ends at $k = 100$ to ensure a fair comparison among the four methodologies. The four methodologies are initialized with a flat start.

Two different scenarios are considered for testing the proposed approach. 1) \emph{Normal operation}: the system's loads change with respect to the aforementioned trend with noisy measurements, and 2) \emph{Sudden load change}: similar to the first scenario but with an arbitrary increment of three times in the load at bus 9 at time step $k = 70$. The increment disappears at time step $k = 71$. After a few tests, we find that having independent ESNs for the voltage magnitude and angle estimation yields the best overall performance. 

% 14-bus test case
\begin{table}[t]
\caption{Voltage angle MAE - Normal operation - IEEE 14-Bus test case} % - Normal Operation
\begin{center}
\begin{tabular}{c||c|c|c|c} 
 \hline
 % &&&&\\
Bus&\multicolumn{4}{c}{Algorithm}\\
 \cline{2-5}
 %Ratio&(\%)&(\%)&(\%)&(\%)\\
Number& ESN & EKF & PF & WLS\\
 %&$k = 10$&$k = 40$&$k = 80$&$k = 130$\\
 \hline
 \hline
  1	& $0.005$	& $0.008$	& $0.063$	& $0.043$	\\
 \hline
 2	& $0.009$	& $0.012$	& $0.071$	& $0.059$	\\
 \hline
 3	& $0.030$	& $0.020$	& $0.104$	& $0.174$	\\ 
 \hline
 4	& $0.012$	& $0.010$	& $0.108$	& $0.141$	\\ 
 \hline 
 5	& $0.010$	& $0.008$	& $0.096$	& $0.124$	\\ 
 \hline 
 6	& $0.027$	& $0.020$	& $0.095$	& $0.264$	\\ 
 \hline 
 7	& $0.022$	& $0.017$	& $0.087$	& $0.228$	\\
 \hline 
8	& $0.022$	& $0.017$	& $0.110$	& $0.254$	\\
\hline 
9	& $0.030$	& $0.020$	& $0.101$	& $0.257$	\\
\hline
10	& $0.031$	& $0.021$	& $0.110$	& $0.281$	\\
\hline
11	& $0.029$	& $0.021$	& $0.105$	& $0.288$	\\
\hline
12	& $0.031$	& $0.022$	& $0.100$	& $0.327$	\\
\hline
13	& $0.032$	& $0.024$	& $0.102$	& $0.329$	\\
\hline
14	& $0.037$	& $0.023$	& $0.097$	& $0.303$	\\ 
 \hline
%  \hline
%  \textbf{Total}	& $0.327$	& $0.243$	& $1.349$	&  $3.072$	\\	
 \end{tabular}
\end{center}
\label{tb_angle_no}
\vspace*{-0.7cm}
\end{table}

\begin{table}[t]
\caption{Voltage magnitude MAE $(\times 10^{-2})$ - Normal operation - IEEE 14-Bus test case} % Normal operation
\begin{center}
\begin{tabular}{c||c|c|c|c} 
 \hline
 % &&&&\\
Bus&\multicolumn{4}{c}{Algorithm}\\
 \cline{2-5}
 %Ratio&(\%)&(\%)&(\%)&(\%)\\
Number& ESN & EKF & PF & WLS\\
 %&$k = 10$&$k = 40$&$k = 80$&$k = 130$\\
\hline
\hline
1	& $0.006$	& $0.013$	& $0.15$	& $0.390$	\\
\hline
2	& $0.058$	& $0.095$	& $0.155$	& $0.447$	\\
\hline
3	& $0.077$	& $0.129$	& $0.165$	& $0.480$	\\ 
\hline
4	& $0.077$	& $0.080$	& $0.181$	& $0.448$	\\ 
\hline 
5	& $0.069$	& $0.071$	& $0.176$	& $0.437$	\\ 
\hline 
6	& $0.040$	& $0.067$	& $0.181$	& $0.482$	\\ 
\hline 
7	& $0.071$	& $0.070$	& $0.184$	& $0.498$	\\
\hline 
8	& $0.041$	& $0.072$	& $0.156$	& $0.472$	\\
\hline 
9	& $0.088$	& $0.069$	& $0.182$	& $0.524$	\\
\hline
10	& $0.087$	& $0.067$	& $0.171$	& $0.554$	\\
\hline
11	& $0.067$	& $0.065$	& $0.150$	& $0.524$	\\
\hline
12	& $0.057$	& $0.068$	& $0.177$	& $0.540$	\\
\hline
13	& $0.065$	& $0.067$	& $0.186$	& $0.510$	\\
\hline
14	& $0.098$	& $0.069$	& $0.167$	& $0.530$	\\ 
\hline
% \hline
% \textbf{Total}	& $0.901$	& $1.002$	& $2.386$	& $6.835$	\\
\end{tabular}
\end{center}
\label{tb_volt_no}
\vspace*{-0.5cm}
\end{table}

\begin{table}[t]
\caption{Voltage angle MAE - Sudden Load Change - IEEE 14-Bus test case}
\begin{center}
\begin{tabular}{c||c|c|c|c} 
 \hline
 % &&&&\\
Bus&\multicolumn{4}{c}{Algorithm}\\
 \cline{2-5}
 %Ratio&(\%)&(\%)&(\%)&(\%)\\
Number& ESN & EKF & PF & WLS\\
 %&$k = 10$&$k = 40$&$k = 80$&$k = 130$\\
 \hline
1   & $0.026$    & $0.036$    & $0.090$    & $0.083$    \\
 \hline
2   & $0.040$    & $0.050$    & $0.106$    & $0.099$    \\
 \hline
3   & $0.098$    & $0.094$    & $0.094$    & $0.202$   \\ 
 \hline
4   & $0.087$    & $0.107$    & $0.096$    & $0.216$   \\ 
 \hline 
5   & $0.075$    & $0.093$    & $0.082$    & $0.197$   \\ 
 \hline 
6   & $0.135$    & $0.178$    & $0.095$    & $0.372$   \\ 
 \hline 
7   & $0.122$    & $0.159$    & $0.103$    & $0.316$   \\
 \hline 
8   & $0.122$    & $0.159$    & $0.102$    & $0.362$   \\
\hline 
9   & $0.140$    & $0.187$    & $0.103$    & $0.358$   \\
\hline
10  & $0.143$    & $0.190$    & $0.097$    & $0.364$   \\
\hline
11  & $0.140$    & $0.185$    & $0.095$    & $0.397$   \\
\hline
12  & $0.144$    & $0.192$    & $0.096$    & $0.438$   \\
\hline
13  & $0.145$    & $0.192$    & $0.108$    & $0.408$   \\
\hline
14  & $0.156$    & $0.210$    & $0.090$    & $0.452$   \\
\hline
% \hline
% \textbf{Total}	& $1.572$	& $2.032$	& $1.356$	& $4.262$	\\
\end{tabular}
\end{center}
\label{tb_angle_slc}
\vspace*{-0.5cm}
\end{table}

\begin{table}[t]
\caption{Voltage magnitude MAE $(\times 10^{-2})$ - Sudden load change - IEEE 14-Bus test case}
\begin{center}
\begin{tabular}{c||c|c|c|c} 
 \hline
 % &&&&\\
Bus&\multicolumn{4}{c}{Algorithm}\\
 \cline{2-5}
 %Ratio&(\%)&(\%)&(\%)&(\%)\\
Number& ESN & EKF & PF & WLS\\
 %&$k = 10$&$k = 40$&$k = 80$&$k = 130$\\
\hline
1   & 0.005     & 0.012     & 0.175     & 0.396     \\
\hline
2   & 0.094     & 0.131     & 0.157     & 0.454     \\
\hline
3   & 0.150     & 0.204     & 0.178     & 0.516     \\ 
\hline
4   & 0.148     & 0.154     & 0.185     & 0.463     \\ 
\hline 
5   & 0.131     & 0.135     & 0.180     & 0.453     \\ 
\hline 
6   & 0.042     & 0.064     & 0.185     & 0.492     \\ 
\hline 
7   & 0.127     & 0.129     & 0.186     & 0.498     \\
\hline 
8   & 0.041     & 0.070     & 0.162     & 0.480     \\
\hline 
9   & 0.164     & 0.152     & 0.198     & 0.530     \\
\hline
10  & 0.160     & 0.146     & 0.178     & 0.560     \\
\hline
11  & 0.109     & 0.109     & 0.151     & 0.525     \\
\hline
12  & 0.080     & 0.088     & 0.183     & 0.545     \\
\hline
13  & 0.099     & 0.098     & 0.190     & 0.520     \\
\hline
14  & 0.181     & 0.156     & 0.174     & 0.610     \\ 
\hline
\end{tabular}
\end{center}
\label{tb_volt_slc}
\vspace*{-0.5cm}
\end{table}

% 300-bus test case

% \begin{table}[t]
% \caption{Voltage angle MAE - Normal operation - 300 Bus} % - Normal Operation
% \begin{center}
% \begin{tabular}{c|c||c|c|c|c} 
%  \hline
%  % &&&&\\
% Type of & Voltage &\multicolumn{4}{c}{Algorithm}\\
%  \cline{3-6}
%  %Ratio&(\%)&(\%)&(\%)&(\%)\\
% Study   &         & ESN & EKF & PF & WLS\\
%  %&$k = 10$&$k = 40$&$k = 80$&$k = 130$\\
%  \hline
%  \hline
%  Normal             &  Angle     & $0.327$ & $0.243$ & $1.349$ & $3.072$ \\
%  \cline{3-6}
%  operation          & Magnitude	& $0.901$ & $1.002$	& $2.381$ & $6.836$ \\
%  \hline
%  Sudden load        &  Angle     & $1.573$ & $2.032$ & $1.357$ & $4.264$ \\
%  \cline{3-6}
%  change             & Magnitude	& $1.531$ & $1.648$	& $2.482$ & $7.042$ \\
%  \hline
% \end{tabular}
% \end{center}
% \label{tb_angle_no}
% \vspace*{-0.7cm}
% \end{table}

\begin{table}[t]
\caption{Accumulated MAE of the buses for the IEEE 300-Bus test case} % - Normal Operation
\begin{center}
\begin{tabular}{c|c||c|c|c|c} 
 \hline
Type of & Voltage &\multicolumn{4}{c}{Algorithm} \\
 \cline{3-6}
Study   &         & ESN & EKF & PF & WLS \\
 \hline
 \hline
 Normal             &  Angle     & $4.635$ & $4.215$ & $9.745$ & $18.36$ \\
 \cline{3-6}
 operation          & Magnitude	& $0.4505$ & $0.501$	& $1.1905$ & $3.418$ \\
 \hline
 Sudden load        &  Angle     & $10.865$ & $13.16$ & $9.785$ & $24.32$ \\
 \cline{3-6}
 change             & Magnitude	& $0.7655$ & $0.824$	& $1.241$ & $3.521$ \\
 \hline
\end{tabular}
\end{center}
\label{results-300}
\vspace*{-0.7cm}
\end{table}

\subsection{Performance Metric}
The methodologies' performance is evaluated using the voltage magnitude and angle mean absolute error (MAE). The MAE is defined as
\begin{subequations}
\begin{align}
e_{V}(i) &= \frac{1}{T}\sum_{k=1}^{T}|V_i(k) - V_{i,o}(k)|, \:\: \forall i=1,\dots,14\\
e_{\theta}(i) &= \frac{1}{T}\sum_{k=1}^{T}|\theta_i(k) - \theta_{i,o}(k)|, \:\: \forall i=1,\dots,14,
\end{align}
\end{subequations}
where $e_{V}(i)$ and $e_{\theta}(i)$ are the mean absolute errors for voltage magnitude and angle at $i$th bus. $V_i(k)$ is the estimated voltage magnitude at time step $k$, $V_{i,o}(k)$ is the true voltage magnitude at time step $k$, $\theta_i(k)$ is the estimated voltage magnitude at time step $k$, $\theta_{i,o}(k)$ is the true voltage magnitude at time step $k$.

\subsection{Hyperparameter Tuning}
The ESN has many hyper-parameters to tune, making it hard to find the best hyper-parameter combination. Typically, one finds the best possible hyper-parameters combination that leads to a non-global optimal solution. Second, the ridge regression solution heavily relies on the true state vector of the previous time step. In machine learning terms, the solution of \eqref{eq:stableSolu} tends to ``over-fit" with the last seen true state vector. Based on these two facts, the ridge regression update of $\bm{W}_k^{\text{out}}$ yields a sub-optimal solution. Hence, the second term of \eqref{eqlms} can be interpreted as a bias that reduces the over-fitting or as an error-correcting term proportional to the deviation of $\bm{W}_k^{\text{out}}$ from the true output. 

In this work, the hyperparameters of the ESN were appropriately tuned for both the 14-bus and 300-bus test cases, according to \cite{Lukosevicius2012}. Similarly, the WLS, EKF, and PF hyperparameters were set based on \cite{Cevallos20181} for both test cases. Regarding the EKF, the initial covariance matrix is set to $10^{-6}$ \cite{Valverde2011}. The parameters $\alpha_k$ and $\beta_k$ are set to 0.8 and 0.5 \cite{LeitedaSilva1983}. The diagonal elements of $\bm{Q}$ and $\bm{R}$ are set to $10^{-6}$. The number of particles in PF was chosen by trial and error, considering the optimal tradeoff between estimation accuracy and processing time.

\subsection{Simulation Results}
% NORMAL OPERATION
Tables \ref{tb_angle_no} and \ref{tb_volt_no} present the comparison of the voltage angle and magnitude absolute errors respectively for normal operating conditions. The prediction step of the EKF has better filtering of the measurement noise, yielding improved estimation performance. This can be seen as low values for the absolute error in the tables. The recurrence property of the ESN and its reservoir's non-linearity can capture the measurements and state dynamics. The ESN is the second-best estimator. It has a performance comparable to the EKF for voltage magnitude estimation. The WLS is expected to be the worst estimator, given that its estimation is based only on the available noisy measurements. 

The methodologies are assessed in the presence of a sudden load change with results shown in Tables \ref{tb_angle_slc} and \ref{tb_volt_slc}. The PF only outperforms the other methods in the voltage angle estimation because the voltage magnitude differences between two consecutive time steps are minor compared to the voltage angle differences. Higher differences yield a better computation of the PF likelihood, which means that the particle weights are different enough, making it possible to select particles with higher weights. On the contrary, more minor differences will apportion similar weights to most particles, becoming difficult to discriminate the particles with the highest weights. The EKF reduces its performance in the sudden load change scenario because the EKF relies on the prediction step based on the previous state estimate, which is very different from the actual system's condition. The ESN is the second-best estimator for the voltage angle and the best for the voltage magnitude. The ESN's echo state property protects the ESN's reservoir against over-excitations driven by the input signal. Similar results can be seen in Table \ref{results-300}, where the four estimation algorithms were evaluated for the IEEE 300-bus test case. Table \ref{results-300} presents the accumulated MAE of all buses for voltage angle and magnitudes for both scenarios, normal operation and sudden load change.

Finally, Tables \ref{tb_time} and \ref{tb_time_300} present the processing time comparison of all the estimation algorithms for both test cases. The ESN is the method with the least processing time in both scenarios and test cases.

\begin{table}[t]
\caption{Comparison of the Computational Time - IEEE 14-Bus test case}
\begin{center}
\begin{tabular}{c||c|c|c|c} 
 \hline
 % &&&&\\
Type of&\multicolumn{4}{c}{Algorithm}\\
 \cline{2-5}
 %Ratio&(\%)&(\%)&(\%)&(\%)\\
Scenario& ESN & EKF & PF & WLS\\
 %&$k = 10$&$k = 40$&$k = 80$&$k = 130$\\
\hline
Normal Operation   & 1.1 s & 14.6 s & 92.2 s & 1.5 s \\
\hline
Sudden Load Change & 1.1 s & 13.4 s & 93.7 s & 1.5 s \\
% \hline
% Bad Data& & & & \\ 
\hline
\end{tabular}
\end{center}
\label{tb_time}
\vspace*{-0.5cm}
\end{table}

\begin{table}[t]
\caption{Comparison of the Computational Time - IEEE 300-Bus test case}
\begin{center}
\begin{tabular}{c||c|c|c|c} 
 \hline
 % &&&&\\
Type of&\multicolumn{4}{c}{Algorithm}\\
 \cline{2-5}
 %Ratio&(\%)&(\%)&(\%)&(\%)\\
Scenario& ESN & EKF & PF & WLS\\
 %&$k = 10$&$k = 40$&$k = 80$&$k = 130$\\
\hline
Normal Operation    & 5.3 s & 27.1 s & 153.5 s & 6.5 s\\
\hline
Sudden Load Change  & 5.2 s & 25.8 s & 158.3 s & 6.7 s\\
% \hline
% Bad Data& & & & \\ 
\hline
\end{tabular}
\end{center}
\label{tb_time_300}
\vspace*{-0.5cm}
\end{table}

\section{Conclusion}\label{sec:conclu}
In this paper, we develop an ESN estimator for electric power systems. Simulated scenarios provide a comprehensive understanding of the advantages of the ESN over industry practices, including the WLS, EKF, and PF. Extensive numerical results show that ESN's model-free method has comparable or even better estimation accuracy than complex model-based methods. More importantly, the ESN reduces estimation time up to $26\%$ and more than $90\%$ compared to the WLS and EKF/PF methods.

\bibliographystyle{IEEEtran}
\bibliography{references}

\end{document}